\documentclass[12pt,reqno]{amsart}
\usepackage{amssymb,amsmath,amsthm,newlfont}
\usepackage{psfrag}
\usepackage[dvips,final]{graphicx}
\numberwithin{equation}{section}

\newcommand{\cI}{{\mathcal{I}}}
\newcommand{\cJ}{{\mathcal{J}}}
\newcommand{\cD}{{\mathcal{D}}}
\newcommand{\cA}{{\mathcal{A}}}
\newcommand{\cL}{{\mathrm{lip}_\alpha}\,}
\newcommand{\C}{\mathbb{C}\,}

\newcommand{\T}{\mathbb{T}\,}
\newcommand{\D}{\mathbb{D}\,}
\newcommand{\N}{\mathbb{N}\,}

\newtheorem{theo}{\bf Theorem}[section]

\newtheorem{rem}[theo]{\bf Remark}
\newtheorem{lem}[theo]{\bf Lemma}

\begin{document}

\title{ {Closed ideals in some algebras
of analytic functions }}

%\date{}

\author[B. Bouya]{Brahim Bouya}

\email{brahimbouya@gmail.com}

\thanks{This work was partially supported by the "Action Integr\'ee
Franco-Marocaine" No.\ MA/03/64.}

\begin{abstract}
{We obtain a complete description of closed ideals of the algebra
$\cD\cap \cL,$ $0<\alpha\leq\frac{1}{2},$ where $\cD$ is the
Dirichlet space and $\cL$ is the algebra of analytic functions
satisfying the Lipschitz condition of order $\alpha.$}
\end{abstract}
\maketitle

\section{Introduction}

The Dirichlet space $\cD$  consists of the complex-valued analytic
functions $f$ on the unit disk $\D$ with finite Dirichlet integral
$$D(f):=\int_{\D}|f'(z)|^2dA(z)<+\infty,$$ where
$dA(z)=\frac{1}{\pi}rdrdt$ denotes the normalized area measure on
$\D.$ Equipped with the pointwise algebraic operations and   the
norm
$$\|f\|^2_{\cD}:= \frac{1}{2\pi}\int_0^{2\pi}|f(e^{it})|^2 dt +
D(f)=\sum_{n=0}^{\infty}(1+n)|\hat{f}(n)|^{2},$$ $\cD$ becomes a
Hilbert space. For $0<\alpha\leq 1,$  let $\cL$ be the algebra of
analytic functions $f$  on $\D$ that are continuous on
$\mathbb{\overline{D}}$ satisfing  the Lipschitz condition of order
$\alpha$ on $\overline{\D}:$
$$|f(z)-f(w)|= o(|z-w|^\alpha)\qquad (\ |z-w|\rightarrow 0).$$
Note that this condition is equivalent to
$$|f'(z)|=o((1-|z|)^{\alpha-1})\qquad   (|z|\rightarrow 1^-).$$
Then, $\cL$ is a Banach algebra when equipped with the norm
$$\|f\|_\alpha:=\|f\|_\infty+ \sup\{(1-|z|)^{1-\alpha}|f'(z)|:
z\in\D\}.$$ Here $\|f\|_\infty:=\sup_{z\in\D}|f(z)|.$ Unlike as for
the case when  $0<\alpha\leq{1}/{2}$, the inclusion $\cL\subset \cD$
always holds provided that   ${1}/{2}< \alpha \leq 1$. In what
follows,  let  $0<\alpha\leq{1}/{2}$ and define $\cA_\alpha:=\cD\cap
\cL$.  It is easy to check that $\cA_\alpha$ is a commutative Banach
algebra when it is endowed with the pointwise algebraic operations
and the norm
$$\|f\|_{\cA_\alpha}:=\|f\|_\alpha + D^{1/2}(f)\qquad (f\in\cA_\alpha).$$
In order to describe the closed ideals in subalgebras of the disc
algebra $A(\D),$ it is natural to make use of Nevanlinna's
factorization theory. For $f\in A(\D)$ there is a canonical
factorization $f=C_f U_f O_f$, where $C_f$ is a constant, $U_f$ a
inner function that is $|U_f|=1$ $a.e$ on $\T$ and $O_f$  the outer
function given by
$$O_f(z)=\exp\Big\{\frac{1}{2\pi}\int_{0}^{2\pi}\frac{e^{i\theta}+z}{e^{i\theta}-z}\log|f(e^{i\theta})|d\theta\Big\}.$$
Denote by $\mathcal{H^{\infty}}(\D)$ the algebra of bounded analytic
functions. Note that $\cA_\alpha$ has the so-called F-property
\cite{Shi2, Car}: if $f\in \cA_\alpha$ and $U$ is an inner function
such that $f/U\in\mathcal{H}^{\infty}(\D)$ then $f/U\in \cA_\alpha$
and $\|f/U\|_{\cA_\alpha}\leq C_{\alpha} \|f\|_{\cA_\alpha}$, where
$C_\alpha$ is independent of $f$. Korenblum \cite{Kor} has described
the closed ideals of the algebra $H^2_1$ of analytic functions $f$
such that $f'\in H^2,$ where $H^2$ is the Hardy space. This result
has been extended to some other Banach algebras of analytic
functions, by Matheson  for $\cL$ \cite{Mat} and by Shamoyan for the
algebra $\lambda^{(n)}_\omega$ of analytic functions $f$ on $\D$
such that
$$|f^{(n)}(\zeta_{1})-
f^{(n)}(\zeta_{2})|=o(\omega(|\zeta_1-\zeta_2|))\ \ as\
|\zeta_1-\zeta_2|\rightarrow 0,$$ where $n$ is a nonnegative integer
and $\omega$ an arbitrary nonnegative nondecreasing subadditive
function on $(0,+\infty)$ \cite{Sha}. Shirokov \cite{Shi1, Shi2} had
given a complete description of closed ideals for  Besov algebras
$AB^{s}_{p,q}$ of analytic functions and particularly for the case
$s>1/2$ and $p=q=2$
$$AB_{2,2}^{s}=\Big\{f\in A(\D)\text{ : }\sum_{n\geq0}|\widehat{f}(n)|^{2}(1+n)^{2s}<\infty\Big\}.$$
Note that the case of $AB_{2,2}^{1/2}=A(\D)\cap\cD$ the problem of
description of closed ideals appears to be much more difficult (see
\cite{HS, EKR}). The purpose of this paper is to describe the
structure of the closed ideals of the Banach algebras $\cA_\alpha.$
More precisely we prove that these ideals are standard in the sense
of the Beurling-Rudin characterization of the closed ideals in the
disc algebra \cite{Hof}:
\begin{theo}
\label{coro} If $\cI$ is closed ideal of $\cA_\alpha$, then
$$\cI=\Big\{f \in\cA_\alpha \text{ : } f_{|E_{_{\cI}}}=0 \text{ and } f/U_{_{\cI}}\in\mathcal{H^{\infty}}(\D)\Big\},$$
where  $E_{_{\cI}}:=\{z\in\T\text{ : }f(z)= 0 ,\ \forall f\in\cI\}$
and $U_{_{\cI}}$ is the greatest common divisor of the inner parts
of the non-zero functions in $\cI.$
\end{theo}
Such characterization of closed ideals can be reduced further to a
problem of approximation of outer functions using the
Beurling--Carleman--Domar resolvent method. Define $d(\xi,E)$ to be
the distance from $\xi\in\T$ to the set $E\subset\T.$ Suppose that
$\cI$ is a closed ideal in $\cA_\alpha$ such that $U_{_{\cI}}=1.$ We
have $Z_{_{\cI}}=E_{_{\cI}},$ where
$$Z_{_{\cI}}:=\{z\in\overline{\D}\text{ : }f(z)= 0 ,\ \forall f\in\cI\}.$$
Next, for $f\in\cA_\alpha$ such that
$$|f(\xi)|\leq C d(\xi,E_{_{\cI}})^{M_\alpha} \qquad (\xi\in\T),$$ where
$M_\alpha$ is a positive constant depending only on $\cA_\alpha,$ we
have  $f\in \cI$  (see section 3 for more precisions). Now, to prove
Theorem \ref{coro} we need Theorem \ref{theo1} below, which states
that every function in $\cA_\alpha\setminus\{0\}$ can be
approximated in $\cA_\alpha$ by functions with boundary zeros of
arbitrary high order.
\begin{theo}\label{theo1}
Let $f$ be a function in $\cA_\alpha\setminus\{0\}$ and let $M>0.$
There exists a sequence of functions $\{g_n\}_{n=1}^\infty\subset
A(\D)$ such
that\\
\begin{enumerate}
\item For all $n\in\N,$ we have $f_n=fg_n\in\cA_\alpha$ and $\displaystyle\lim_{n\rightarrow\infty} \|f_n - f\|_{\cA_\alpha}=0.$
\\
\item  $|g_n(\xi)|\leq C_n d^{M}(\xi,E_f)\qquad(\xi\in\T),$\\
\end{enumerate}
where $E_f:=\{\xi\in\T\text{ : } f(\xi)=0\}.$
\end{theo}
To prove this Theorem, we give a refinement of the classical
Korenblum approximation theory \cite{Kor, Mat, Sha, Shi1, Shi2}.

\section{ main result on approximation of functions in $\cA_\alpha$}

We begin by fixing some notations. Let $f\in \cA_\alpha$ and let
$\{\gamma_{n}:=(a_n,b_n) \}_{n\geq0}$ be  the countable collection
of the (disjoint open) arcs of $\T\setminus E_f.$ Without loss of
the generality, we can suppose that the arc lengths of $\gamma_n$
are less than $1/2$. In what follows, we denote by $\Gamma$ the
union of a family of arcs $\gamma_{n}.$ Define
$$f_{_{\Gamma}}(z):=\exp\Bigm\{\frac{1}{2\pi}\int_{\Gamma}\frac{e^{i\theta}+z}{e^{i\theta}-z}\log|f(e^{i\theta})|d\theta
\Bigm\}.$$ The difficult part in the proof of Theorem \ref{theo1} is
to establish the following

\begin{theo}\label{theo2}
Let $f\in\cA_\alpha\setminus\{0\}$ be an outer function such that
$\|f\|_{\cA_\alpha}\leq 1$ and let $N\geq1$ and $\rho>1.$ Then we
have $f^\rho f^{N}_{_{\Gamma}}\in \cA_\alpha$ and
\begin{eqnarray}\label{gol}\sup_{ \Gamma}\| f^\rho f^{N}_{_{\Gamma}}\|_{\cA_\alpha}\leq
C_{N,\rho},\end{eqnarray} where $C_{N,\rho}$ is a positive constant
independent of $\Gamma.$
\end{theo}

\begin{rem}
\label{rem} \rm{ For a set $S\subset A(\D),$ we denote by $co(S)$
the convex hull of $S$ consisting of the intersection of all convex
sets that contain $S.$ Set $\Gamma_n=\cup_{m\geq n}\gamma_m$ and let
$f$ be as in the Theorem \ref{theo2}. It is clear that the sequence
$(f^\rho f^{N}_{_{\Gamma_n}})_n$ converges uniformly on compact
subsets of $\D$ to $f^\rho.$ We use \eqref{gol} to deduce, by the
Hilbertian structure of $\cD$, that there is a sequence $h_n\in
co(\{f^\rho f^{N}_{_{\Gamma_m}}\}_{m=n}^{\infty})$ converging to
$f^\rho$ in ${\cD}$. Also, by \cite[section4]{Mat}, we obtain that
$h_n$ converges to $f^\rho$ in $\cL,$ for sufficiently large $N$ (in
fact, we can prove that this result remains true for every
$N\geq1$). Therefore $\|h_n - f^\rho\|_{\cA_\alpha}\to 0,$ as $n\to
\infty.$}
\end{rem}

Define $\cJ(F)$ to be the closed ideal of all functions in
$\cA_\alpha$ that vanish on $F\subset\overline{\D}.$ In the proof of
Theorem \ref{theo1}, we need the following classical lemma, see for
instance \cite[Lemma 4]{Mat} and \cite[Lemma 24]{Kor}.

\begin{lem}\label{finitezeros}
Let $f\in\cA_\alpha$ and $E'$ be a finite subset of  $\T$ such that
$f_{|E'}=0$. Let $M>0$ be given. For every $\varepsilon>0$ there is
an outer function $F$ in $\cJ(E')$ such that
\begin{enumerate}
\item $  \|Ff - f\|_{\cA_\alpha}\leq\varepsilon,$\\
\item $ |F(\xi)|\leq C d^{M}(\xi,E') \qquad (\xi\in\T).$
\end{enumerate}
\end{lem}
\vspace{0.3cm} \textit{\bf Proof of Theorem \ref{theo1} :} Now, we
can deduce the proof of Theorem \ref{theo1} by using Theorem
\ref{theo2} and Lemma \ref{finitezeros}. Indeed, let $f$ be a
function in $\cA_\alpha\setminus\{0\}$ such that
$\|f\|_{\cA_\alpha}\leq 1$ and let $\epsilon>0.$  For $m\geq 1$ we
have
$$(fO_{f}^{\frac{1}{m}}-f)'=(O_{f}^{\frac{1}{m}}-1)f'+\frac{1}{m}U_{f}O_{f}^{\frac{1}{m}}O_{f}^{'}.$$
The F-property of $\cA_\alpha$ implies that $O_f \in \cA_\alpha$.
Then, there exists $\eta_0\in\N$ such that
$$\|fO_{f}^{\frac{1}{m}}-f\|_{\cA_\alpha}<\epsilon/3\qquad (m\geq \eta_0).$$
Set $\Gamma_n=\cup_{p\geq n}\gamma_p$ and $N\geq M/\alpha$ for a
given $M>0.$ By Remark \ref{rem} applied to $O_{f}$ (with
$\rho=1+\frac{1}{m}$), there is a sequence $k_{n,m}\in
co(\{f_{_{_{\Gamma_p}}}^{N}\}_{p=n}^\infty)$ such that
$$\|O_{f}^{1+\frac{1}{m}}k_{n,m}-O_{f}^{1+\frac{1}{m}}\|_{\cA_\alpha}<\frac{1}{m}\qquad (n\in\N,\ m\geq1).$$
It is clear that
$$\|O_{f}^{\frac{1}{m}}f_{_{_{\Gamma_n}}}^{N}-O_{f}^{\frac{1}{m}}\|_{\infty}\longrightarrow 0 \qquad(n\longrightarrow+\infty).$$
Then for every $m\geq1$ we get
$$\|O_{f}^{\frac{1}{m}}k_{n,m}-O_{f}^{\frac{1}{m}}\|_{\infty}\longrightarrow 0 \qquad(n\longrightarrow+\infty).$$
So, there is a sequence $k_{m}\in
co(\{f_{_{_{\Gamma_p}}}^{N}\}_{p=m}^\infty)$ such that
$$\left\{
  \begin{array}{ll}
   \|O_{f}^{1+\frac{1}{m}}k_{m}-O_{f}^{1+\frac{1}{m}}\|_{\cA_\alpha}\leq\frac{1}{m}  & (m\geq1),\\
\|O_{f}^{\frac{1}{m}}k_{m}-O_{f}^{\frac{1}{m}}\|_{\infty}\leq
\frac{1}{m} & (m\geq1).
  \end{array}
\right.$$
We have
$$(fO_{f}^{\frac{1}{m}}k_m-fO_{f}^{\frac{1}{m}})'=(f'-U_fO_f^{'})
(O_{f}^{\frac{1}{m}}k_m-O_{f}^{\frac{1}{m}})+U_f(O_{f}^{1+\frac{1}{m}}k_m-O_{f}^{1+\frac{1}{m}})'.$$
Since $\|O_f\|_{\cA_\alpha}\leq C_\alpha\|f\|_{\cA_\alpha}\leq
C_\alpha,$ we obtain
\begin{eqnarray*}
&&\|fO_{f}^{\frac{1}{m}}k_m-fO_{f}^{\frac{1}{m}}\|_{\cA_\alpha}\\
&=&\|fO_{f}^{\frac{1}{m}}k_m-fO_{f}^{\frac{1}{m}}\|_{\infty}+\sup_{z\in\D}\{(1-|z|)^{1-\alpha}|(fO_{f}^{\frac{1}{m}}k_m-fO_{f}^{\frac{1}{m}})'(z)|\}
\\&&+D^{1/2}(fO_{f}^{\frac{1}{m}}k_m-fO_{f}^{\frac{1}{m}})\\
&\leq& \|fO_{f}^{\frac{1}{m}}k_m-fO_{f}^{\frac{1}{m}}\|_{\infty}
+C_\alpha\|f\|_{\alpha}
\|O_{f}^{\frac{1}{m}}k_m-O_{f}^{\frac{1}{m}}\|_{\infty}
\\&&+\sup_{z\in\D}\{(1-|z|)^{1-\alpha}|(O_{f}^{1+\frac{1}{m}}k_m-O_{f}^{1+\frac{1}{m}})'(z)|\}
\\&&+C
\|O_{f}^{\frac{1}{m}}k_m-O_{f}^{\frac{1}{m}}\|_{\infty}D^{1/2}(f)
+CD^{1/2}(O_{f}^{1+\frac{1}{m}}k_m-O_{f}^{1+\frac{1}{m}})
\\&\leq&
C_\alpha\|O_{f}^{\frac{1}{m}}k_m-O_{f}^{\frac{1}{m}}\|_{\infty}+
C\|O_{f}^{1+\frac{1}{m}}k_m-O_{f}^{1+\frac{1}{m}}\|_{\cA_\alpha}
\\&\leq& \frac{C_\alpha}{m}.
\end{eqnarray*}
Then, fix $\eta_1\geq \eta_0$ such that
$$\|fO_{f}^{\frac{1}{m}}k_m-fO_{f}^{\frac{1}{m}}\|_{\cA_\alpha}<\epsilon/3\qquad (m\geq \eta_1).$$
We have $k_m=\sum\limits_{i\leq j_m} c_i f_{_{_{\Gamma_{i}}}}^{N},$
where $\sum\limits_{i\leq j_m}c_i=1.$ Set $E'_m=\cup_{i<j_m}\partial
\gamma_i.$ Using Lemma \ref{finitezeros}, we obtain an outer
function $F_{m}\in \cJ(E^{'}_{m})$  such that $ |F_{m}(\zeta)|\leq
C_{m} d^{M}(\zeta,E^{'}_{m})$ for $\zeta\in \T$ and
$$\|f O_{f}^{\frac{1}{m}}k_m F_{m}-f O_{f}^{\frac{1}{m}}k_m\|_{\cA_\alpha}<\frac{1}{m}\qquad (m\geq 1).$$
Then fix $\eta_2\geq \eta_1$ such that
$$\|f O_{f}^{\frac{1}{m}}k_m F_{m}-f O_{f}^{\frac{1}{m}}k_m\|_{\cA_\alpha}<\epsilon/3\qquad (m\geq \eta_2).$$
Consequently we obtain $$\|fO_{f}^{\frac{1}{m}}k_m F_{m}-
f\|_{\cA_\alpha}<\epsilon\qquad (m\geq \eta_2).$$ It is not hard to
see that
$$|O_{f}^{\frac{1}{m}}k_m F_m(\xi)|\leq|k_m F_m(\xi)|\leq C_{m}d^{M}(\xi,E_f)\qquad (\xi\in \T).$$ Therefore
$g_m =O_{f}^{\frac{1}{m}}k_m F_{m}$ is the desired sequence, which
completes the proof of Theorem \ref{theo1}.

\section{Beurling--Carleman--Domar resolvent method }

Since $\cA_\alpha\subset \cL$, then for all $f\in \cA_\alpha,$ $E_f$
satisfies the Carleson condition
$$\int_\T\log\frac{1}{d(e^{it},E_f)}dt<+\infty.$$
For $f\in\cA_\alpha,$ we denote by $B_f$ the Blashke product with
zeros $Z_{f}\setminus E_{f},$ where
$Z_{f}:=\{z\in\overline{\D}\text{ : } f(z)=0\}.$ We begin with
following lemma

\begin{lem}\label{lemdiv}
Let $\cI$ be a closed ideal of $\cA_\alpha.$ Define $B_{_\cI}$ to be
the Blashke product with zeros $Z_{_\cI}\setminus E_{_\cI}.$ There
is a function $f\in\cI$ such that $B_f=B_{_\cI}.$
\end{lem}

\begin{proof}
Let $g\in\cI$ and let $B_n$ be the Blashke product with zeros
$Z_{g}\cap\D_{n},$ where $\D_{n}:=\{z\in\D\text{ :
}|z|<\frac{n-1}{n},\ n\in\N\}.$ Set $g_n=g/K_n,$ where $K_n=B_n/I_n$
and $I_n$ is the Blashke product with zeros $Z_{\cI}\cap\D_{n}.$ We
have $g_n\in\cI$ for every $n.$ Indeed, fix $n\in\N.$  It is
permissible to assume that $Z_{K_n}$ consists of a single point, say
$Z_{K_n}=\{w\}.$ Let $\pi:\cA_\alpha\to \cA_\alpha/\cI$ be the
canonical quotient map. First suppose $w\notin Z_{_\cI},$ then
$\pi(K_n)$ is invertible in $\cA_\alpha/\cI.$ It follows that
$\pi(g_n)=\pi(g)\pi^{-1}(K_n)=0,$ hence $g_n\in\cI.$ If $w\in
Z_{_\cI},$ we consider the following ideal
$\cJ_w:=\{f\in\cA_\alpha\text{ : }fI_n\in\cI\}.$  It is clear that
$\cJ_w$ is closed. Since $w\notin Z_{_{\cJ_w}},$ it follows that
$K_n$ is invertible in the quotient algebra $\cA_{\alpha}/\cJ_w$ and
so $g/(I_nK_n)\in\cJ_w.$ Hence $g_n\in\cI.$ \\
It is clear that $g_n$ converges uniformly on compact subsets of
$\D$ to $f=(g/B_g)B_\cI$ and we have $B_f=B_{_{\cI}}.$ In the sequel
we prove that $f\in\cI.$ If we obtain
$$|\big(g_n\big)'(z)|\leq o\Big(\frac{1}{(1-r)^{1-\alpha}}\Big)\qquad(z\in\D),$$
uniformly with respect to $n,$ we can deduce by using \cite[Lemma
1]{Mat} that
$\lim\limits_{n\rightarrow+\infty}\|g_n-f\|_{\alpha}=0.$ Indeed, by
the Cauchy integral formula
\begin{eqnarray*}
\big(g_n\big)'(z)&=&\frac{1}{2\pi i}\int_\T
\frac{g(\zeta)\overline{K_n(\zeta)}}{(\zeta-z)^2}d\zeta\\&=&\frac{1}{2\pi
i}\int_\T
\frac{(g(\zeta)-g(z/|z|))\overline{K_n(\zeta)}}{(\zeta-z)^2}d\zeta\qquad(z\in\D).
\end{eqnarray*}
Then, for $z=re^{i\theta}\in\D$
\begin{eqnarray*}
|\big(g_n\big)'(z)|&\leq& \frac{\|K_n\|_{\infty}}{2\pi}\int_\T
\frac{|g(\zeta)-g(z/|z|)|}{|\zeta-z|^2}|d\zeta|\\
&=&\frac{1}{2\pi}\int_{-\pi}^{\pi}\frac{|g(e^{i(t+\theta)})-g(e^{i\theta})|}
{1-2r\cos t+r^2}dt.
\end{eqnarray*}
For all $\varepsilon>0,$ there is $\eta>0$ such that if
$|t|\leq\eta,$ we have
$$|g(e^{i(t+\theta)})-g(e^{i\theta})|\leq \varepsilon|t|^\alpha\qquad(\theta\in[-\pi,+\pi]).$$
Then
\begin{eqnarray*}
&&\int_{-\pi}^{\pi}\frac{|g(e^{i(t+\theta)})-g(e^{i\theta})|}
{1-2r\cos t+r^2}dt\\&\leq&\varepsilon
\int_{|t|\leq\eta}\frac{|t|^{\alpha}}{(1-r)^2+4rt^2/\pi^2}dt+
\|g\|_{\alpha}\int_{|t|\geq\eta}\frac{|t|^{\alpha}}{(1-r)^2+4rt^2/\pi^2}dt
\\&\leq&
\frac{\varepsilon}{r^{\frac{1+\alpha}{2}}(1-r)^{1-\alpha}}\int_{0}^{+\infty}\frac{u^{\alpha}}{1+(2u/\pi)^2}du
\\&&\hspace{3.5cm}+\frac{\|g\|_{\alpha}}{r^{\frac{1+\alpha}{2}}(1-r)^{1-\alpha}}\int_{|u|\geq\frac{\eta\sqrt{r}}{1-r}}\frac{u^{\alpha}}{1+(2u/\pi)^2}du
\\&\leq& \varepsilon
O\Big(\frac{1}{(1-r)^{1-\alpha}}\Big)+\|g\|_{\alpha}o\Big(\frac{1}{(1-r)^{1-\alpha}}\Big).
\end{eqnarray*}
We obtain
\begin{equation}\label{unifmaj}
\int_{-\pi}^{\pi}\frac{|g(e^{i(t+\theta)})-g(e^{i\theta})|}
{1-2r\cos
t+r^2}dt\leq\|g\|_{\alpha}o\Big(\frac{1}{(1-r)^{1-\alpha}}\Big).
\end{equation}
Consequently
$$|\big(g_n\big)'(z)|\leq\|g\|_{\alpha}o\Big(\frac{1}{(1-r)^{1-\alpha}}\Big)\qquad(z\in\D).$$
By the F-property of $\cA_\alpha,$ we have $\|g_n\|\leq
C_\alpha\|g\|_{\cA_\alpha}.$ Using the Hilbertian structure of
$\cD,$ we deduce that there is a sequence $h_n\in
co(\{g_k\}_{k=n}^{\infty})$ converging to $f$ in ${\cD}$. It is
clear that $h_n\in\cI$ and
$\lim\limits_{n\rightarrow+\infty}\|h_n-f\|_{\alpha}=0.$ Then
$\lim\limits_{n\rightarrow+\infty}\|h_n-f\|_{\cA_\alpha}=0.$ Thus
$f\in\cI.$ This completes the proof of the lemma.
\end{proof}

 As a consequence of Theorem
\ref{theo1}, we can prove Theorem \ref{coro} and deduce that each
closed ideal of $\cA_\alpha$ is standard. For the sake of
completeness, we sketch here the proof.

\vspace{0.3cm} \textit{\bf Proof of Theorem \ref{coro} :} Define
$\gamma$ on $\D$ by $\gamma(z)=z$ and let $\pi:\cA_\alpha\to
\cA_\alpha/\cI$ be the canonical quotient map. Also, let
$f\in\cJ(E_{_{\cI}})$ be such that
$f/U_{_{\cI}}\in\mathcal{H}^{\infty}(\D) $ and $({f_n})_n$ be the
sequence in Theorem \ref{theo1} associated to $f$ with $M\geq 3$.
More exactly, we have $f_n=fg_n,$ where $|g_n(\xi)|\leq
d^3(\xi,E_f)\leq d^3(\xi,E_{_{\cI}}).$ Define
$$\mathrm{L_\lambda}(f)(z):=\left\{
\begin{array}{ll}
\displaystyle \frac{f(z)-f(\lambda)}{z-\lambda}\ & \mbox{if}\ z\neq\lambda,\\
\displaystyle f'(\lambda)\ & \mbox{if}\ z=\lambda.\\
\end{array}\right.$$
Then
\begin{equation}\label{lambda}\pi(f)(\pi(\gamma)-\lambda)^{-1}=f(\lambda)(\pi(\gamma)-\lambda)^{-1}+\pi\big(\mathrm{L_\lambda}(f)\big).
\end{equation}
It is clear that $(\pi(\gamma)-\lambda)^{-1}$ is an analytic
function on $\C\setminus Z_{_{\cI}}.$ Note that the multiplicity of
the pole $z_0\in Z_{_{\cI}}\cap\D$ of $(\pi(\gamma)-\lambda)^{-1}$
is equal to the multiplicity of the zero $z_0$ of $U_{_{\cI}}.$
Since $U_{_{\cI}}$ divides $f,$ then according to \eqref{lambda} we
can deduce that $\pi(f)(\pi(\gamma)-\lambda)^{-1}$ is an analytic
function on $\C\setminus E_{_{\cI}}.$ Let $|\lambda|>1,$ we have
\begin{eqnarray}\nonumber
\|\pi(f)(\pi(\gamma)-\lambda)^{-1}\|_{\cA_\alpha}&\leq&
\|f\|_{\cA_\alpha}\sum_{n=0}^{\infty}\|\gamma^{n}\|_{\cA_\alpha}|\lambda|^{-n-1}\\\label{der1}&\leq&
\|f\|_{\cA_\alpha}\frac{C}{(|\lambda|-1)^{3/2}}.
\end{eqnarray}
By Lemma \ref{lemdiv}, there is $g\in\cI$ such that
$B_g=B_{_{\cI}}.$ Let $k=f(g/B_g).$ Then $k=(f/B_{_{\cI}})g\in\cI$
and for $|\lambda|<1,$ we have
$$k(\lambda)(\pi(\gamma)-\lambda)^{-1}=-\pi\big(\mathrm{L_\lambda}(k)\big).$$
Therefore
\begin{eqnarray} \nonumber
\|\pi(f)(\pi(\gamma)-\lambda)^{-1}\|_{\cA_\alpha}&\leq&
|f(\lambda)|\|(\pi(\gamma)-\lambda)^{-1}\|_{\cA_\alpha}+\|\mathrm{L_\lambda}(f)\|_{\cA_\alpha}
\\&\leq& \nonumber\frac{\|\mathrm{L_\lambda}(k)\|_{\cA_\alpha}}{|g/B_g|(\lambda)}+\|\mathrm{L_\lambda}(f)\|_{\cA_\alpha}
\\&\leq& \nonumber \frac{C(f,k)}{(1-|\lambda|)|g/B_g|(\lambda)}
\\&\leq& \label{der2}C(f,k)e^{\frac{C}{1-|\lambda|}}\quad(|\lambda|<1).
\end{eqnarray}
We use \cite[Lemmas 5.8 and 5.9]{TW} to deduce
$$\|\pi(f)(\pi(\gamma)-\xi)^{-1}\|\leq \frac{C(f,k)}{d(\xi,E_{_{\cI}})^3}\qquad (1\leq|\xi|\leq 2,\ \xi\notin E_{_{\cI}}).$$
Then, we obtain
$$\xi\mapsto|(g_n)(\xi)|\|\pi(f)(\pi(\gamma)-\xi)^{-1}\|\in
L^\infty(\T).$$ With a simple calculation as in \cite[Lemma
2.4]{ESZ}, we can deduce that
$$\pi({f_n})=\frac{1}{2\pi i}\int_\T (g_n)(\xi)\pi(f)(\pi(\gamma)-\xi)^{-1} d\xi.$$
Denote  $\cI^{\infty}_{U_{_{\cI}}}(E_{_{\cI}}):=\{h\in A(\D)\text{ :
} h_{|E_{_{\cI}}}=0\ \mbox{and}\ h/U_{_{\cI}}\in A(\D)\}.$ From
\cite[p. 81]{Hof}, we know that
$\cI^{\infty}_{U_{_{\cI}}}(E_{_{\cI}})$ has an approximate identity
$(e_m)_{m\geq1}\in\cI^{\infty}_{U_{_{\cI}}}(E_{_{\cI}})$ such that
$\| e_m\|_\infty\leq 1.$  ${\cI}$ is dense in
$\cI^{\infty}_{U_{_{\cI}}}(E_{_{\cI}})$ with respect to the sup norm
$\|\hspace{0.1cm}.\hspace{0.1cm}\|_{\infty}$, so there exists
$(u_m)_{m\geq 1}\in \cI$ with $\|u_m\|_\infty\leq 1$ and
$\displaystyle\lim_{m\to\infty}u_m(\xi)=1$ for $\xi\in\T\setminus
E_{_{\cI}}$. Therefore $\pi(f_n)=\pi(f_n-f_nu_m)\to 0$ as $m\to
\infty$. Then $f_n\in \cI$ and $f\in \cI.$

\section{Proof of Theorem \ref{theo2}}

The proof of Theorem \ref{theo2} is based on a series of lemmas. In
what follows, $C_\rho$ will denote a positive number that depends
only on $\rho,$ not necessarily the same at each occurrence.
 For an open subset  $\Delta$
of $\D$, we put
$$\|f'\|_{L^2(\Delta)}^{2}:=\int_{\Delta}|f'(z)|^{2}dA(z).$$
We begin with the following key lemma

\begin{lem}\label{local} Let $f\in \cA_\alpha$ be such that $\|f\|_{\cA_\alpha}\leq1$ and let
$\rho>1$ be given. Then
$$\int_{\gamma}\frac{|f(e^{it})|^{2\rho}}{d(e^{it})}dt\leq
C_{\rho}\|f'\|^2_{L^2(\Delta_\gamma)},$$ where $a,b\in E_f,$
$\gamma=(a,b)\subset\T\setminus E_f,$ $d(z):=\min\{|z-a|,|z-b|\}$
and  $\Delta_\gamma:=\{z\in\D:\ z/|z|\in \gamma\}.$
\end{lem}

\begin{proof} Let $e^{it}\in\gamma$ and define $z_t:=(1-d(e^{it}))e^{it}.$
Since $|\gamma|<1/2,$ we obtain $|z_t|>1/2.$ We have
\begin{equation}\label{loc1}
|f(e^{it})|^{2\rho}\leq2^{2\rho-1}\Bigm(|f(e^{it})-f(z_t)|^{2\rho}+|f(z_t)|^{2\rho}\Bigm).
\end{equation}
By  H\"older's inequality combined with the fact that
$\|f\|_{\infty}\leq\|f\|_{\cA_\alpha}\leq1,$ we get
\begin{eqnarray*}
|f(e^{it})-f(z_t)|^{2\rho}&=&|f(e^{it})-f(z_t)|^{2\rho-2}|f(e^{it})-f(z_t)|^{2}\\
&\leq& 2^{2\rho-2}(1-|z_t|)\int_{|z_t|}^{1}|f'(re^{it})|^{2}dr\\
&\leq& 2^{2\rho-1}d(e^{it})\int_{0}^{1}|f'(re^{it})|^{2}rdr.
\end{eqnarray*}
Hence
\begin{eqnarray}\nonumber
\int_{\gamma}\frac{|f(e^{it})- f(z_t)|^{2\rho}}{d(e^{it})}dt &\leq&
2^{2\rho-1}\int_{\gamma}\int_{0}^{1}|f'(re^{it})|^{2}rdrdt
\\\label{loc2} &\leq&
2^{2\rho-1}\pi\|f'\|^2_{L^2(\Delta_\gamma)}.
\end{eqnarray}
Since $d(e^{it})\leq1/2$, we obtain $\frac{d(e^{it})}{\sqrt{2}}\leq
d(z_t)\leq \sqrt{2}d(e^{it})$. Put $d(z_t)=|z_t-\xi|$ and note that
either $\xi=a$ or $\xi=b.$ Let
$$z_{t}(u)=(1-u)z_t+u\xi  \qquad (0\leq u\leq1).$$
With a simple calculation, we can prove that for all
$e^{it}\in\gamma$ and for all $u,$ $0\leq u\leq1,$ we have
$$|z_{t}(u)-w|>\frac{1}{2}(1-u)d(e^{it})\qquad (w\in\partial\Delta_\gamma),$$
where $\partial\Delta_\gamma$ is the boundary of $\Delta_\gamma.$
Then $\D_{t,u}:=\{z\in\D :
|z-z_{t}(u)|\leq\frac{1}{2}(1-u)d(e^{it})\}\subset\Delta_\gamma,$
for all $e^{it}\in\gamma$ and for all $u,$  $0\leq u\leq1.$ Since
$|f'(z)|$ is subharmonic on $\D$, it follows that
\begin{eqnarray*}
|f'(z_{t}(u))|&\leq&
\frac{4}{\pi(1-u)^2d^2(e^{it})}\int_{\D_{t,u}}|f'(z)|dA(z)\\&\leq&\frac{2}{\pi^{1/2}(1-u)d(e^{it})}
\|f'\|_{L^2(\Delta_\gamma)}.
\end{eqnarray*}
Set $\varepsilon_{\rho}=2\alpha(\rho-1).$ We have
\begin{eqnarray*}
|f^{\rho}(z_t)|^{2}&=&|f^{\rho}(z_t)-f^{\rho}(\xi)|^{2}\\&=&
\rho^2|z_t-\xi|^2|\int_{0}^{1}f^{\rho-1}(z_{t}(u))f'(z_{t}(u))du|^2\\
&\leq& C_\rho
d^2(e^{it})\Bigm(\int_{0}^{1}|z_{t}(u)-\xi|^{\frac{\varepsilon_\rho}{2}}|f'(z_{t}(u))|du\Bigm)^2\\
&\leq& C_\rho
d^{\varepsilon_\rho}(e^{it})\Bigm(\int_{0}^{1}\frac{1}{(1-u)^{1-\frac{\varepsilon_\rho}{2}}}du\Bigm)^2
\|f'\|^2_{L^2(\Delta_\gamma)}\\
&\leq& C_\rho d^{\varepsilon_\rho}(e^{it})
\|f'\|^2_{L^2(\Delta_\gamma)}.
\end{eqnarray*}
Hence
\begin{eqnarray}
\int_{\gamma}\frac{|f(z_t)|^{2\rho}}{d(e^{it})}dt\leq
C_\rho\|f'\|^2_{L^2(\Delta_\gamma)}. \label{loc3}
\end{eqnarray}
Therefore the result follows from  \eqref{loc1}, \eqref{loc2} and
\eqref{loc3}.
\end{proof}

In the sequel, we denote by $f$ an outer function in $\cA_\alpha$
such that $\|f\|_{\cA_\alpha}\leq1$ and we fix a constant $\rho,$
$1<\rho\leq 2.$ By \cite[Theorem B]{Mat}, we have $f^\rho
f^{N}_{_{\Gamma}}\in \cL$ and $\|f^\rho
f^{N}_{_{\Gamma}}\|_{\cL}\leq C_{N,\rho}$. To prove Theorem
\ref{theo2} we need to estimate the integral
$\int_\D|(f^{\rho}f_{_{\Gamma}}^{N})'|^{2}dA(z).$  Define
\begin{equation}
g_{_{\Gamma}}(z):=\frac{1}{\pi}\int\limits_{\Gamma}\frac{e^{i\theta}}{(e^{i\theta}-z)^{2}}\log
|f( e^{i\theta})| d\theta.
\end{equation}
Clearly we have $f'=f(g_{_{\Gamma}}+g_{_{\T\setminus\Gamma}})$ and
$(f_{_{\Gamma}}^{N})'=Nf_{_{\Gamma}}^{N}g_{_{\Gamma}}$, so
\begin{eqnarray}\label{eq2}f^{\rho}(f_{_{\Gamma}}^{N})'&=&Nf^{\rho}f_{_{\Gamma}}^{N}g_{_{\Gamma}}\\\label{eq3}&=&
f^{\rho-1}Nf'f_{_{\Gamma}}^{N}-Nf^{\rho}f_{_{\Gamma}}^{N}g_{_{\T\setminus\Gamma}}.
\end{eqnarray}
Since  $\|f\|_{\infty}\leq1$, it is obvious that
$\|f_{_{\Gamma}}^{N}\|_{\infty}\leq1$ and
$\|f^{\rho-1}\|_{\infty}\leq1.$ Hence, by \eqref{eq2} we get
\begin{equation}\label{objet}\int_\D|(f^{\rho}f_{_{\Gamma}}^{N})'|^{2}dA(z)\leq
\rho^2 +
N^{2}\int_\D|f^{\rho}(f_{_{\Gamma}})'|^{2}dA(z).\end{equation} We
fix $\gamma=(a,b)\subset\T\setminus E_f$ such that $f(a)=f(b)=0.$
Our purpose in what follows is to estimate the integral
\begin{equation}\int_{\Delta_\gamma}|f^{\rho}(f_{_{\Gamma}})'|^{2}dA(z)\end{equation}
which we can rewrite as
$$\int_{\Delta_\gamma}|f^{\rho}(f_{_{\Gamma}})'|^{2}
dA(z)=\int_{\Delta_\gamma^{1}}+\int_{\Delta_\gamma^{2}},$$
where \begin{eqnarray*} \Delta_\gamma^{1}&:=&\Bigm\{z\in \Delta_\gamma : \ d(z)< 2(1-|z|)\Bigm\}\\
\Delta_\gamma^{2}&:=&\Bigm\{z\in \Delta_\gamma : \ d(z)\geq
2(1-|z|)\Bigm\}.
\end{eqnarray*}

\subsection{The integral on the region $\Delta_\gamma^{1}$.}

We begin with the following lemma.

\begin{lem}\label{lem2}
$$\int_{\Delta_\gamma}\frac{|f(z)-f(z/|z|)|^{2\rho}}{(1-|z|)^2}dA(z)\leq \frac{1}{2\alpha(\rho-1)}
\|f'\|^2_{L^2(\Delta_\gamma)}.$$
\end{lem}

\begin{proof}  Let $z=re^{it}\in \Delta_\gamma$ and put
$\varepsilon_{\rho}=2\alpha(\rho-1).$  We have
\begin{eqnarray*}
r|f(re^{it})-f(e^{it})|^{2\rho}&=&
r|f(re^{it})-f(e^{it})|^{2\rho-2}|f(re^{it})-f(e^{it})|^{2}\\
&\leq&
r(1-r)^{1+\varepsilon_{\rho}}\int_{r}^{1}|f'(se^{it})|^{2}ds\\
&\leq& (1-r)^{1+\varepsilon_{\rho}}\int_{0}^{1}|f'(se^{it})|^{2}sds.
\end{eqnarray*}
Therefore
\begin{eqnarray*}
&&\int_{\Delta_\gamma}\frac{|f(z)-f(z/|z|)|^{2\rho}}{(1-|z|)^2}dA(z)\\
&=&
\int_{0}^{1}\Bigm(\int_\gamma|f(re^{it})-f(e^{it})|^{2\rho}\frac{rdt}{\pi}\Bigm)\frac{dr}{(1-r)^2}\\
&\leq&
\|f'\|^2_{L^2(\Delta_\gamma)}\int_{0}^{1}\frac{1}{(1-r)^{1-\varepsilon_{\rho}}}dr.
\end{eqnarray*}
This completes the proof.
\end{proof}

Now, we can state the following result.

\begin{lem}\label{D1}
$$\int_{\Delta_\gamma^{1}}|f(z)|^{2\rho}|f'_{_{\Gamma}}(z)|^{2}dA(z)\leq
C_{\rho}\|f'\|^2_{L^2(\Delta_\gamma)}.$$
\end{lem}

\begin{proof}  By Cauchy's estimate, it follows that
$|f'_{_{\Gamma}}(re^{it})|\leq \frac{1}{1-r}.$ Using  Lemma
\ref{lem2}, we get
\begin{eqnarray}\nonumber
&&\int_{\Delta_\gamma^{1}}|f(z)|^{2\rho}|f'_{_{\Gamma}}(z)|^{2}dA(z)\\\nonumber&\leq&
\int_{\Delta_\gamma^{1}}\frac{|f(z)|^{2\rho}}{(1-|z|)^2}dA(z)\\
\label{G11}&\leq& C_{\rho}\|f'\|^2_{L^2(\Delta_\gamma)}+
2^{2\rho-1}\int_{\Delta_\gamma^{1}}\frac{|f(z/|z|)|^{2\rho}}{(1-|z|)^2}dA(z).
\end{eqnarray}
Using Lemma \ref{local}, we obtain
\begin{eqnarray}\nonumber
\int_{\Delta_\gamma^{1}}\frac{|f(z/|z|)|^{2\rho}}{(1-|z|)^2}dA(z)&=&\frac{1}{\pi}
\int_{\Delta_\gamma^{1}}\frac{|f(e^{it})|^{2\rho}}{(1-r)^2}rdrdt
\\\nonumber
&\leq& \frac{C}{\pi}\int_{\gamma}
\frac{|f(e^{it})|^{2\rho}}{d(e^{it})}dt
\\\label{G12}
&\leq& C_{\rho}\|f'\|^2_{L^2(\Delta_\gamma)}.
\end{eqnarray}
The result of our lemma follows by combining the estimates
(\ref{G11}) and (\ref{G12}).
\end{proof}

\subsection{The integral on the region $\Delta_\gamma^{2}$}

In this subsection, we estimate the integral
$\int_{\Delta_\gamma^{2}}|f(z)|^{2\rho}|f'_{_{\Gamma}}(z)|^{2}dA(z).$
Before this, we make some remarks. For  $z\in\D$ define
$$a_{\gamma}(z):=\left\{
\begin{array}{ll}
\displaystyle\frac{1}{2\pi}\int_\Gamma\frac{-\log|f(e^{i\theta})|}{|e^{i\theta}-z|^{2}}d\theta \quad&\mbox{if}\ \gamma\nsubseteq\Gamma,\\
\displaystyle\frac{1}{2\pi}\int_{\T\setminus\Gamma}\frac{-\log|f(e^{i\theta})|}{|e^{i\theta}-z|^{2}}d\theta \quad &\mbox{if}\ \gamma\subseteq\Gamma.\\
\end{array}\right.
$$
Using the equation (\ref{eq2}), it is easy to see that
\begin{equation}\label{naturel}|f(z)^{\rho}f'_{_{\Gamma}}(z)|^{2}\leq
4\Bigm|f(z)^{\rho}\frac{1}{2\pi}\int_\Gamma\frac{-\log|f(e^{i\theta})|}{|e^{i\theta}-z|^{2}}d\theta\Bigm|^{2}.
\end{equation}
Using  the equation (\ref{eq3}), it is clear that
\begin{equation}\label{naturel2}|f(z)^{\rho}f'_{_{\Gamma}}(z)|^{2}\leq
2|f'(z)|^{2} +
8\Bigm|f(z)^{\rho}\frac{1}{2\pi}\int_{\T\setminus\Gamma}\frac{-\log|f(e^{i\theta})|}{|e^{i\theta}-z|^{2}}d\theta\Bigm|^{2}.
\end{equation}

Then

\begin{multline}\label{simple}
\int_{\Delta_\gamma^{2}}|f(z)|^{2\rho}|f'_{_{\Gamma}}(z)|^{2}dA(z)\leq
2\|f'\|^2_{L^2(\Delta_\gamma)}\\
+8\int_{\Delta_\gamma^{2}}|f(z)|^{2\rho}a_{\gamma}^2(z)dA(z).
\end{multline}

Since $\log|f|\in L^1(\T),$ we have
\begin{equation}
a_{\gamma}(z)\leq\frac{C}{d^2(z)} \qquad (z\in\Delta_\gamma).
\label{evid}
\end{equation}
Given such inequality, it is not easy to estimate immediately the
integral of the function $|f(z)|^{2\rho}a^{2}_{\gamma}(z)$ on the
whole $\Delta_\gamma^{2}.$ In what follows, we give a partition of
$\Delta_\gamma^{2}$ into three parts so that one can  estimate  the
integral $\int|f(z)|^{2\rho}a^{2}_{\gamma}(z)dA(z)$ on each part.
Let $z\in\Delta_\gamma^{2},$ three situations are possible :
\begin{eqnarray} \label{gauche} &a_{\gamma}(z)\leq
\displaystyle8\frac{|\log(d(z))|}{d(z)},\\& \label{milieu}
\displaystyle8\frac{|\log(d(z))|}{d(z)}< a_{\gamma}(z)<
\displaystyle8\frac{|\log(d(z))|}{1-r},\\& \label{droite}
\displaystyle8\frac{|\log(d(z))|}{1-r}\leq a_{\gamma}(z).
\end{eqnarray}
We can now divide $\Delta_\gamma^{2}$ into the following three parts
\begin{eqnarray*} \Delta_\gamma^{21}&:=&\Bigm\{z\in
\Delta_\gamma^{2} : z\ \mbox{satisfying} \ (\ref{gauche}) \Bigm\},\\
\Delta_\gamma^{22}&:=&\Bigm\{z\in \Delta_\gamma^{2} : z\
\mbox{satisfying} \
(\ref{milieu})\Bigm\},\\
\Delta_\gamma^{23}&:=&\Bigm\{z\in \Delta_\gamma^{2} : z\
\mbox{satisfying} \ (\ref{droite})\Bigm\}.
\end{eqnarray*}

\subsubsection{\bf The integral on the regions $\Delta_\gamma^{21}$ and $\Delta_\gamma^{23}$}

In this case we begin by the following

\begin{lem}\label{D21}
$$\int_{\Delta_\gamma^{21}}|f(z)|^{2\rho}a^{2}_{\gamma}(z)dA(z)\leq C_{\rho}\|f'\|^2_{L^2(\Delta_\gamma)}.$$
\end{lem}

\begin{proof} Using Lemma \ref{lem2}, we get
\begin{eqnarray*}
&&\int_{\Delta_\gamma^{21}}|f(z)|^{2\rho}a^{2}_{\gamma}(z)dA(z)\\
 &\leq& 2^\rho\int_{\Delta_\gamma^{21}}
|f(z)|^{\rho-1}|f(z)-f(z/|z|)|^{\rho+1}a^{2}_{\gamma}(z)dA(z)\\
&&\qquad\qquad\qquad\quad +2^\rho\int_{\Delta_\gamma^{21}}|f(z)|^{\rho-1}|f(z/|z|)|^{\rho+1}a^{2}_{\gamma}(z)dA(z)\\
&\leq&
C_{\rho}\int_{\Delta_\gamma}\frac{|f(z)-f(z/|z|)|^{\rho+1}}{(1-|z|)^2}dA(z)
+C_{\rho}
\int_{\Delta_\gamma^{21}}\frac{|f(e^{it})|^{\rho+1}}{d^{2}(e^{it})}rdrdt
\\
&\leq& C_{\rho}\|f'\|^2_{L^2(\Delta_\gamma)}
+C_\rho\int_{\Delta_\gamma^{21}}\frac{|f(e^{it})|^{\rho+1}}{d^{2}(e^{it})}drdt=I_{2,1}.
\end{eqnarray*}
Let $e^{it}\in\gamma$ and  denote by  $\zeta_t$ the point of
$\partial \Delta_\gamma^{2}\cap \D$  such that
$\zeta_t/|\zeta_t|=e^{it}$. We have
$$|e^{it}-\zeta_t|=1-|\zeta_t|=\frac{d(\zeta_t)}{2}\leq d(e^{it}).$$
Then
\begin{eqnarray*}
\int_{\Delta_\gamma^{21}}\frac{|f(e^{it})|^{\rho+1}}{d^{2}(e^{it})}drdt
&\leq&
\int_{\Delta_\gamma^{2}}\frac{|f(e^{it})|^{\rho+1}}{d^{2}(e^{it})}drdt\\
&=& \int_{\gamma}\frac{|f(e^{it})|^{\rho+1}}{d^2(e^{it})}
\int_{|\zeta_t|}^{1}dr  dt \\&\leq&
\int_{\gamma}\frac{|f(e^{it})|^{\rho+1}}{d(e^{it})}dt.
\end{eqnarray*}
Using Lemma \ref{local}, we get
\begin{equation*}
I_{2,1}\leq C_{\rho}\|f'\|^2_{L^2(\Delta_\gamma)}.
\end{equation*}
This proves the result.
\end{proof}

\begin{lem}\label{D23}
$$\int_{\Delta_\gamma^{23}}|f(z)|^{2\rho}a^{2}_{\gamma}(z)dA(z)\leq C A(\Delta_\gamma),$$
where $A(\Delta_\gamma)$ is the area measure of $\Delta_\gamma.$
\end{lem}

\begin{proof}Set
$$\Lambda_\gamma:=\left\{
  \begin{array}{ll}
    \Gamma \quad& \hbox{for}\ \gamma\nsubseteq\Gamma,\\
    \T\setminus\Gamma \quad & \hbox{for}\ \gamma\subseteq\Gamma.
  \end{array}
\right.$$ Let $z\in\Delta_\gamma^{23}.$ We have
\begin{eqnarray*} |f(z)|&=&\exp\Bigm\{\frac{1}{2\pi}\int_{0}^{2\pi}\frac{1-r^{2}
}{|e^{i\theta}- z|^{2}}\log|f(e^{i\theta})|d\theta\Bigm\}\\
&\leq&
\exp\Bigm\{\frac{1}{2\pi}\int_{\Lambda_\gamma}\frac{1-r}{|e^{i\theta}-z|^{2}}\log|f(e^{i\theta})|d\theta\Bigm\}
\\  &=& \exp\Bigm\{-(1-r)a_{\gamma}(z)\Bigm\}
\\  &\leq& d^{8}(z).
\end{eqnarray*}
Using \eqref{evid}, we obtain the result.
\end{proof}

\subsubsection{\bf The integral on the region $\Delta_\gamma^{22}$}

Here, we will give an estimate of the following integral
$$\int_{\Delta_\gamma^{22}}|f(z)|^{2\rho}a^{2}_{\gamma}(z)dA(z).$$
Before doing this, we begin with some lemmas. The next one is
essential for what follows. Note that a similar result is used by
different authors: Korenblum \cite{Kor}, Matheson \cite{Mat},
Shamoyan \cite{Sha} and Shirokov \cite{Shi1,Shi2}.

\begin{lem}\label{lem6} Let $z\in\Delta_\gamma^{22}$ and let
$\mu_{z}=1-\frac{8|\log(d(z))|}{a_{\gamma}(z)}.$ Then
\begin{eqnarray}\label{mu}
|f(\mu_{z} z)|\leq d^{2}(z).
\end{eqnarray}
\end{lem}

\begin{proof}  Let $z\in\Delta_\gamma$ and let $\mu<1.$ We have
\begin{eqnarray*} |f(\mu
z)|&=&\exp\Bigm\{\frac{1}{2\pi}\int_{0}^{2\pi}\frac{1-(\mu r)^{2}
}{|e^{i\theta}-\mu z|^{2}}\log|f(e^{i\theta})|d\theta\Bigm\}\\
&\leq& \exp\Bigm\{\frac{1}{2\pi}\int_{\Lambda_\gamma}\frac{1-(\mu
r)^{2}}{|e^{i\theta}-\mu z|^{2}}\log|f(e^{i\theta})|d\theta\Bigm\}
\\  &\leq& \exp\Bigm\{-(1-\mu
r)\inf_{\theta\in\Lambda_\gamma}\Bigm|\frac{e^{i\theta}-z}{e^{i\theta}-\mu
z}\Bigm|^{2}a_{\gamma}(z)\Bigm\}.
\end{eqnarray*}
For $z\in\Delta_\gamma^{22},$ it is clear that $1-\mu_{z}\leq
d(z)\leq |e^{i\theta}-z|$ for all $e^{i\theta}\in\Lambda_\gamma.$
Then
$$\inf\limits_{\theta\in\Lambda_\gamma}\Bigm|\frac{e^{i\theta}-z}{e^{i\theta}-\mu_{z}
z}\Bigm|\geq \frac{1}{2}\qquad(z\in\Delta_\gamma^{22}).$$ Thus
$$ |f(\mu_{z} z)|\leq \exp\Bigm\{-\frac{1-\mu_{z} }{4}a_{\gamma}(z)\Bigm\}\qquad(z\in\Delta_\gamma^{22}).$$
Then, we have
 $$|f(\mu_{z} z)|\leq\exp\Bigm\{-\frac{1}{4}(1-\mu_{z} )a_{\gamma}(z)\Bigm\}= d^{2}(z)\qquad(z\in\Delta_\gamma^{22}),$$
which yields (\ref{mu}).
\end{proof}

For $r<1$, define
$$\gamma_{r}:=\{z\in\D:\ |z|=r\ \mbox{and}\
z/|z|\in\gamma\}.$$ Without loss of generality, we can suppose that
$d(z)\leq\frac{1}{2},$ $z\in \Delta_\gamma^{2}.$ We need the
following

\begin{lem}\label{lem7}Let $r<1.$ Then
$$
\int_{\gamma_{r}\cap\Delta_\gamma^{22}}|f(re^{it})-f(\mu_{re^{it}}
re^{it})|^{2\rho}a^{2}_{\gamma}(re^{it})rdt
\leq\frac{C_\rho}{(1-r)^{1-\varepsilon_\rho}}\|f'\|^2_{L^2(\Delta_{\gamma})},
$$ where $\varepsilon_\rho=\alpha(\rho-1).$
\end{lem}

\begin{proof}Let $re^{it}\in\Delta_\gamma^{22}.$ Then
\begin{multline*}
|f(re^{it})-f(\mu_{re^{it}} re^{it})|^{\rho-1}[(1-\mu_{re^{it}}
)a_{\gamma}(re^{it})]^{2}\\ \leq 64 (1-\mu_{re^{it}}
)^{\varepsilon_\rho}\log^2(d(re^{it}))\leq C_\rho.
\end{multline*}
It is clear that $1-r\leq 1-\mu_{re^{it}}\leq d(re^{it})\leq
\frac{1}{2}$ and so $\frac{1}{2}\leq\mu_{re^{it}}\leq r.$ We have
\begin{eqnarray*}
&&\int_{\gamma_{r}\cap\Delta_\gamma^{22}}|f(re^{it})-f(\mu_{re^{it}}
re^{it})|^{2\rho}a^{2}_{\gamma}(re^{it})rdt \\
&\leq& C_\rho\int_{\gamma_{r}\cap
\Delta_\gamma^{22}}\frac{|f(re^{it})-f(\mu_{re^{it}}
re^{it})|^{\rho+1}}{(1-\mu_{re^{it}})^{2}
}rdt\\
&\leq& \frac{C_{\rho}}{(1-r
)^{1-\varepsilon_\rho}}\int_{\gamma_{r}\cap
\Delta_\gamma^{22}}\frac{|f(re^{it})-f(\mu_{re^{it}}
re^{it})|^{2}}{1-\mu_{re^{it}}}rdt\\
&\leq& \frac{C_{\rho}}{(1-r
)^{1-\varepsilon_\rho}}\int_{\gamma_{r}\cap \Delta_\gamma^{22}}
\Bigm(\int_{\mu_{re^{it}} r}^{r}|f'(se^{it})|^{2}ds\Bigm)rdt\\
&\leq&\frac{C_{\rho}}{(1-r )^{1-\varepsilon_\rho}}
\int_{S_{r}}|f'(se^{it})|^{2}sdsdt\\&\leq&
\frac{C_{\rho}}{(1-r)^{1-\varepsilon_\rho}}\int_{S_{r}}|f'(w)|^{2}dA(w),
\end{eqnarray*}
where $$S_{r}:=\Bigm\{w\in\D :\ 0\leq |w|\leq r\ \mbox{and}\ \
\frac{w}{|w|}\in\gamma \Bigm\}.$$ The proof is therefore completed.
\end{proof}

The last result that we need before giving the proof of Theorem
\ref{theo2} is the following one.

\begin{lem}\label{D22}
$$\int_{\Delta_\gamma^{22}}|f(z)|^{2\rho}a^{2}_{\gamma}(z)dA(z)\leq C_\rho\|f'\|^2_{L^2(\Delta_{\gamma})}+ C A(\Delta_\gamma).$$
\end{lem}

\begin{proof} Using  \eqref{evid} and Lemmas \ref{lem6} and \ref{lem7},  we find that
\begin{eqnarray*}&&
\int_{\Delta_\gamma^{22}}|f(z)|^{2\rho}a^{2}_{\gamma}(z)dA(z)\\&=&
\frac{1}{\pi}\int_0^1\Bigm(\int_{ \gamma_{r}\cap\Delta_\gamma^{22}}|f(re^{it})|^{2\rho}a_{\gamma}^2(re^{it})rdt\Bigm)dr\\
&\leq& C A(\Delta_\gamma)\\&&+
2^{2\rho-1}\int_0^1\Bigm(\int_{\gamma_{r}\cap\Delta_\gamma^{22}}|f(re^{it})-f(\mu_{re^{it}}
re^{it})|^{2\rho}a_{\gamma}^2(re^{it})rdt\Bigm)dr\\
&\leq& C A(\Delta_\gamma)+ C_\rho\|f'\|^2_{L^2(\Delta_\gamma)}.
\end{eqnarray*}
This completes the proof of the lemma.
\end{proof}

\subsubsection{\bf Conclusion}

Now, according to \eqref{simple} and Lemmas \ref{D21}, \ref{D23} and
\ref{D22}, we obtain
\begin{eqnarray*}
\int_{ \Delta_\gamma^{2}}|f(z)|^{2\rho}|f'_{_{\Gamma}}(z)|^{2}dA(z)
&\leq& 2\|f'\|^2_{L^2(\Delta_\gamma)}+
8\int_{\Delta_\gamma^{2}}|f(z)|^{2\rho}a^{2}_{\gamma}(z)dA(z)
\\ &\leq& C_{\rho}\|f'\|^2_{L^2(\Delta_\gamma)}+ C A(\Delta_\gamma).
\end{eqnarray*}
Combining this with Lemma \ref{D1}, we deduce that
\begin{eqnarray*}
\int_{ \Delta_\gamma}|f(z)|^{2\rho}|f'_{_{\Gamma}}(z)|^{2}dA(z)\leq
C_{\rho}\|f'\|^2_{L^2(\Delta_\gamma)}+ C A(\Delta_\gamma).
\end{eqnarray*}
Hence
\begin{eqnarray*}
\int_{\D}|f(z)|^{2\rho}|f'_{_{\Gamma}}(z)|^{2}dA(z)&=&\sum_{n=1}^{\infty}\int_{\Delta_{\gamma_{n}}}|f(z)|^{2\rho}|f'_{_{\Gamma}}(z)|^{2}dA(z)\\
&\leq&
C_{\rho}\sum_{n=1}^{\infty}\|f'\|^2_{L^2(\Delta_{\gamma_{n}})}+ C
\sum_{n=1}^{\infty}A(\Delta_{\gamma_{n}})\\ &\leq& C_{\rho}.
\end{eqnarray*}
This completes the proof of Theorem \ref{theo2}.

 \vspace{1em}

\textsc{Acknowledgements.} I wish to thank Professors A. Borichev,
O. El Fallah and K. Kellay for the interest which they carried to
this work.

\end{document}